\documentclass[12pt, reqno]{amsart}
\allowdisplaybreaks
\begin{document}
\setcounter{page}{1}

\title[\hfilneg \hfil Gruss type inequality  ]
{Some Gruss type Inequalities using Pompeiu's mean value theorem on Conformable Fractional Calculus}

\author[Deepak B. Pachpatte\hfil \hfilneg]
{Deepak B. Pachpatte}

\address{Deepak B. Pachpatte \newline
 Department of Mathematics,
 Dr. Babasaheb Ambedkar Marathwada University, Aurangabad,
 Maharashtra 431004, India}
\email{pachpatte@gmail.com}


\begin{abstract}
           The main objective of this paper is to obtain some Gruss Like Inequalities using Pompeiu's mean value theorem on Conformable Fractional Calculus.
	
 \end{abstract}
\maketitle

\section{\textbf{Introduction}}
   \paragraph{} In 1935 G. Gruss \cite{Gru} proved the following integral inequality:
	\[
\begin{array}{l}
 \left| {\frac{1}{{b - a}}\int\limits_a^b {f\left( x \right)g\left( x \right)} dx - } \right|\left( {\frac{1}{{b - a}}\int\limits_a^b {f\left( x \right)} dx} \right)\left( {\frac{1}{{b - a}}\int\limits_a^b {g\left( x \right)} dx} \right) \\ 
  \le \frac{1}{4}\left( {P - p} \right)\left( {Q - q} \right), \\ 
 \end{array}
\]
for all $x \in [a,b]$. The constant $\frac{1}{4}$ is the best possible.
\paragraph{}In 1946 D. Pompeiu \cite{Pop} proved the variant of following mean value theorem: 
   For every real valued function f differentiable on an interval $[a,b]$ not containing $0$ and for all pairs $x_1  \ne x_2 $ in $[a,b]$, there exist a point $\xi $ between $x_1$ and $x_2$ such that
	\[
\frac{{x_1 f\left( {x_2 } \right) - x_2 f\left( {x_1 } \right)}}{{x_1  - x_2 }} = f\left( \xi  \right) - \xi f'\left( \xi  \right).
\]
	\paragraph{}Recently \cite{Abd,Ben, Kha} have introduced  the conformable fractional calculus. Since then  many authors have studied the various inequalities and other problems on conformable fractional calculus see \cite{Iyi,Jia,Ust}. The main aim of our paper is to study Some Gruss type Inequalities using Pompeiu's mean value theorem on Conformable Fractional Calculus.
	
	\section{\textbf{Preliminaries}}
	Fractional Inequalities are very effective tool for studying wide range of problems in various branches of mathematics. This has attracted the  attention of large number of mathematicians.
	These type of  Inequalities arise in variety of applications and its study is interesting. 
	Motivated by the above results in this paper we obtain the Gruss type inequality on Conformable Fractional Calculus.
	\paragraph{}Now we give some basic definitions required for proving our results. The Conformable fractional derivative is given as 
	\paragraph{\textbf{Definition 2.1}} \cite{Erd} For a given function $f:\left[ {0,\infty } \right) \to \mathbb{R}$ the conformable fractional derivative of order $0 < \alpha  \le 1$ of $f$ at $t>0$ is defined by
	\[
D_\alpha  \left( f \right)\left( t \right) = \mathop {\lim }\limits_{\varepsilon  \to 0} \frac{{f\left( {t + \varepsilon t^{1 - \alpha } } \right) - f\left( t \right)}}{\varepsilon }.
\]
\paragraph{}The Conformable fractional Integration is given as
		\paragraph{\textbf{Definition 2.2}} \cite{Erd} Let $\alpha  \in (0,1]$ and $0 \le a < b$. A function $f:[a,b] \to \mathbb{R}$ is $\alpha$-fractional integrable on $[a,b]$ if the integral
	\[
\int\limits_a^b {f\left( x \right)d_\alpha  x: = } \int\limits_a^b {f\left( x \right)x^{\alpha  - 1} dx} 
\]
exists and is finite.
		\paragraph{}In order to prove our results we require following result  in \cite{Erd} which is  Pompeiu's Mean value theorem for Conformable Fractional Differentiable function
	\paragraph{\textbf{Theorem 2.1}}\cite{Erd} Let $\alpha  \in (0,1]$, $f:[a,b] \subseteq \mathbb{R} \to \mathbb{R}$ be an $\alpha$-fractional differentiable mapping on $(a,b)$, with $0<a<b$ and for all pairs $x_1  \ne x_2$ in $[a,b]$, there exist a point $\xi$ in $(x_1,x_2)$ such that the following equality holds:
	\[
\frac{{x_1^\alpha  f\left( {x_2 } \right) - x_2^\alpha  f\left( {x_1 } \right)}}{{\frac{{x_1^\alpha  }}{\alpha } - \frac{{x_2^\alpha  }}{\alpha }}} = \alpha f\left( \xi  \right) - \xi ^{2 - \alpha } D_\alpha  \left( f \right)\left( \xi  \right).
\]
 
\section{\textbf{MAIN RESULTS}}
We denote following Notations for our convenience 
\[
\begin{array}{l}
 K[p,q] = \int\limits_a^b {p\left( x \right)} q\left( x \right)d_\alpha  x - \frac{1}{{b^{2\alpha }  - a^{2\alpha } }}\left[ {\left( {\int\limits_a^b {p\left( x \right)d_\alpha  x} } \right)} \right.\left( {\int\limits_a^b {x^\alpha  q\left( x \right)d_\alpha  x} } \right) \\ 
 \left. { + \left( {\int\limits_a^b {q\left( x \right)d_\alpha  x} } \right)\left( {\int\limits_a^b {x^\alpha  p\left( x \right)d_\alpha  x} } \right)} \right],  
 \end{array}
\tag{3.1}\]
\[
\begin{array}{l}
 H[p,q] = \int\limits_a^b {p\left( x \right)} q\left( x \right)d_\alpha  x - \frac{{3\alpha }}{{b^{3\alpha }  - a^{3\alpha } }}\left( {\int\limits_a^b {x^\alpha  q\left( x \right)d_\alpha  x} } \right)\left( {\int\limits_a^b {t^\alpha  p\left( t \right)d_\alpha  t} } \right) \\ 
  - \frac{{3\alpha }}{{b^{3\alpha }  - a^{3\alpha } }}\left( {\int\limits_a^b {x^\alpha  p\left( x \right)d_\alpha  x} } \right)\left( {\int\limits_a^b {t^\alpha  q\left( t \right)d_\alpha  t} } \right) + \int\limits_a^b {p\left( t \right)q\left( t \right)} d_\alpha  t . 
 \tag{3.2}\end{array}
\]

Now we give some Gruss type inequality for conformable fractional differentiable function using Pompeiu type inequality as follows
\paragraph{\textbf{Theorem 3.1}}
Let $f,g:[a,b] \to \mathbb{R}$ be continuous on $[a,b]$ and $\alpha$-fractional differential mapping on $(a,b)$ with $0<a<b$. Then
\begin{align*}
\left| {K\left[ {f,g} \right]} \right|
&\le \left\| {f - w_1 D_\alpha  f} \right\|_\infty  \int\limits_a^b {\left| {g\left( x \right)} \right|} \left| {\frac{1}{{2^\alpha  }} - \frac{x}{{a^\alpha   + b^\alpha  }}} \right|d_\alpha  x \\
&+ \left\| {g - w_1 D_\alpha  g} \right\|_\infty  \int\limits_a^b {\left| {f\left( x \right)} \right|} \left| {\frac{1}{{2^\alpha  }} - \frac{x}{{a^\alpha   + b^\alpha  }}} \right|d_\alpha  x ,
\tag{3.3}
\end{align*}
where 
$w_1 \left( t \right) = \frac{{\xi ^{2 - \alpha } }}{\alpha }$, $w_2 \left( t \right) = \frac{{\eta ^{2 - \alpha } }}{\alpha }$,$t \in [a,b]$
and we have
\[
\left\| {f - w_1 D_\alpha  f} \right\|_\infty   = \mathop {\sup }\limits_{\xi  \in \left( {a,b} \right)} \left| {f\left( \xi  \right) - w_1 D_\alpha  \left( f \right)\left( \xi  \right)} \right| < \infty,  
\tag{3.4}\]
and
\[
 \left\| {g - w_1 D_\alpha  g} \right\|_\infty   = \mathop {\sup }\limits_{\xi  \in \left( {a,b} \right)} \left| {g\left( \eta  \right) - w_1 D_\alpha  \left( g \right)\left( \eta  \right)} \right| < \infty.   
\tag{3.5}\]

\paragraph{\textbf{Proof.}}
From Theorem 2.1 we have 
\[
t^\alpha  f\left( x \right) - x^\alpha  f\left( t \right) = \left[ {\frac{{t^\alpha  }}{\alpha } - \frac{{x^\alpha  }}{\alpha }} \right]\left( {\alpha f\left( \xi  \right) - \xi ^{2 - \alpha } D_\alpha  \left( f \right)\left( \xi  \right)} \right),
\tag{3.6}\]
\[
t^\alpha  g\left( x \right) - x^\alpha  g\left( t \right) = \left[ {\frac{{t^\alpha  }}{\alpha } - \frac{{x^\alpha  }}{\alpha }} \right]\left( {\alpha g\left( \eta  \right) - \eta ^{2 - \alpha } D_\alpha  \left( g \right)\left( \eta  \right)} \right).
\tag{3.7}\]

Multiplying $(3.6)$ and $(3.7)$ by $g(x)$ and $f(x)$ respectively, we have
\[
t^\alpha  f\left( x \right)g\left( x \right) - x^\alpha  f\left( t \right)g\left( x \right) = \left[ {\frac{{t^\alpha  }}{\alpha } - \frac{{x^\alpha  }}{\alpha }} \right]\left( {\alpha f\left( \xi  \right) - \xi ^{2 - \alpha } D_\alpha  \left( f \right)\left( \xi  \right)} \right)g\left( x \right),
\tag{3.8}\]
\[
t^\alpha  f\left( x \right)g\left( x \right) - x^\alpha  g\left( t \right)f\left( x \right) = \left[ {\frac{{t^\alpha  }}{\alpha } - \frac{{x^\alpha  }}{\alpha }} \right]\left( {\alpha g\left( \eta  \right) - \eta ^{2 - \alpha } D_\alpha  \left( g \right)\left( \eta  \right)} \right)f\left( x \right).
\tag{3.9}\]
Adding $(3.8)$ and $(3.9)$ we have
\begin{align*}
&2t^\alpha  f\left( x \right)g\left( x \right) - x^\alpha  f\left( t \right)g\left( x \right) - x^\alpha  f\left( x \right)g\left( t \right) \\ 
&= \left( {\frac{{t^\alpha  }}{\alpha } - \frac{{x^\alpha  }}{\alpha }} \right)\left( {\alpha f\left( \xi  \right) - \xi ^{2 - \alpha } D_\alpha  f\left( \xi  \right)} \right)g\left( x \right) \\ 
&+ \left( {\frac{{t^\alpha  }}{\alpha } - \frac{{x^\alpha  }}{\alpha }} \right)\left( {\alpha g\left( \eta  \right) - \eta ^{2 - \alpha } D_\alpha  g\left( \eta  \right)} \right)f\left( x \right).
\tag{3.10}
\end{align*}
Integrating $(3.10)$ on both sides with respect to $t$ from $a$ to $b$ for conformable fractional integral we have
\begin{align*}
&\left( {b^{2\alpha }  - a^{2\alpha } } \right)f\left( x \right)g\left( x \right) - x^\alpha  g\left( x \right)\int\limits_a^b {f\left( t \right)d_\alpha  } t - x^\alpha  f\left( x \right)\int\limits_a^b {g\left( t \right)d_\alpha  } t \\ 
& = \left( {\alpha f\left( \xi  \right) - \xi ^{2 - \alpha } D_\alpha  f\left( \xi  \right)} \right)g\left( x \right)\left( {\frac{{b^{\alpha  + 1}  - a^{\alpha  + 1} }}{{\alpha  + 1}}} \right) \\
&- \left( {\alpha f\left( \xi  \right) - \xi ^{2 - \alpha } D_\alpha  f\left( \xi  \right)} \right)g\left( x \right)\frac{{x^\alpha  }}{\alpha }\left( {b^\alpha   - a^\alpha  } \right) \\ 
&+ \left( {\alpha g\left( \eta  \right) - \eta ^{2 - \alpha } D_\alpha  g\left( \eta  \right)} \right)f\left( x \right)\left( {\frac{{b^{\alpha  + 1}  - a^{\alpha  + 1} }}{{\alpha  + 1}}} \right) \\
&- \left( {\alpha g\left( \eta  \right) - \eta ^{2 - \alpha } D_\alpha  g\left( \eta  \right)} \right)f\left( x \right)\frac{{x^\alpha  }}{\alpha }\left( {b^\alpha   - a^\alpha  } \right).
\tag{3.10}
\end{align*}
Integrating both sides of $(3.10)$ with respect to $x$ over $[a,b]$ we have

\begin{align*}
&\left( {b^{2\alpha }  - a^{2\alpha } } \right)\int\limits_a^b {f\left( x \right)g\left( x \right)} d_\alpha  x - \left( {\int\limits_a^b {f\left( t \right)d_\alpha  } t} \right)\left( {\int\limits_a^b {x^\alpha  g\left( x \right)d_\alpha  x} } \right) \\
&- \left( {\int\limits_a^b {g\left( t \right)d_\alpha  } t} \right)\left( {\int\limits_a^b {x^\alpha  g\left( x \right)d_\alpha  x} } \right) \\ 
& = \left( {\alpha f\left( \xi  \right) - \xi ^{2 - \alpha } D_\alpha  f\left( \xi  \right)} \right) \\ 
& \left\{ {\frac{{\left( {b^{2\alpha }  - a^{2\alpha } } \right)}}{{2\alpha }}\int\limits_a^b {g\left( x \right)d_\alpha  x}  - \frac{{\left( {b^\alpha   - a^\alpha  } \right)}}{\alpha }\int\limits_a^b {x^\alpha  g\left( x \right)d_\alpha  x} } \right\} \\
&  + \left( {\alpha g\left( \eta  \right) - \eta ^{2 - \alpha } D_\alpha  g\left( \eta  \right)} \right) \\
& \left\{ {\frac{{\left( {b^{2\alpha }  - a^{2\alpha } } \right)}}{{2\alpha }}\int\limits_a^b {f\left( x \right)d_\alpha  x}  - \frac{{\left( {b^\alpha   - a^\alpha  } \right)}}{\alpha }\int\limits_a^b {x^\alpha  f\left( x \right)d_\alpha  x} } \right\}.
\tag{3.11}
\end{align*}
Now from $(3.1)$ we have
\begin{align*}
K\left( {p,q} \right)
&  = \left( {\alpha f\left( \xi  \right) - \xi ^{2 - \alpha } D_\alpha  f\left( \xi  \right)} \right)\int\limits_a^b {g\left( x \right)} \left\{ {\frac{1}{{2^\alpha  }} - \frac{x}{{a^\alpha   + b^\alpha  }}} \right\}d_\alpha  x \\
&+ \left( {\alpha g\left( \eta  \right) - \eta ^{2 - \alpha } D_\alpha  g\left( \eta  \right)} \right)\int\limits_a^b {f\left( x \right)} \left\{ {\frac{1}{{2^\alpha  }} - \frac{x}{{a^\alpha   + b^\alpha  }}} \right\}d_\alpha  x .
\tag{3.12}
\end{align*}
From the properties of modulus and $(3.12)$ we have
\begin{align*}
 \left| {K\left( {p,q} \right)} \right| 
&\le \left\| {f - w_1 D_\alpha  f} \right\|_\infty  \int\limits_a^b {\left| {g\left( x \right)} \right|} \left| {\frac{1}{{2^\alpha  }} - \frac{x}{{a^\alpha   + b^\alpha  }}} \right|d_\alpha  x \\
&+ \left\| {g - w_2 D_\alpha  g} \right\|_\infty  \int\limits_a^b {\left| {f\left( x \right)} \right|} \left| {\frac{1}{{2^\alpha  }} - \frac{x}{{a^\alpha   + b^\alpha  }}} \right|d_\alpha  x.
\tag{3.13}
\end{align*}
which is required inequality.
\paragraph{\textbf{Theorem 3.2}}
Let $f,g:[a,b] \to \mathbb{R}$ be continuous on $[a,b]$ and $\alpha$-fractional differential mapping on $(a,b)$ with $[a,b]$ not containing $0$. Then
\begin{align*}
 \left| {H\left( {f,g} \right)} \right| 
&\le \left| {\alpha f\left( \xi  \right) - \xi ^{2 - \alpha } D_\alpha  f\left( \xi  \right)} \right| \times  \\
& \left| {\alpha g\left( \eta  \right) - \eta ^{2 - \alpha } D_\alpha  g\left( \eta  \right)} \right|\left( {\frac{{3\alpha }}{{b^{3\alpha }  - a^{3\alpha } }}} \right) \times  \\
&\left\{ {\frac{{b^{3\alpha }  - a^{3\alpha } }}{{3\alpha }}\left( {b - a} \right) - \frac{{\left( {b^{2\alpha }  - a^{2\alpha } } \right)^2 }}{{2\alpha ^2 }} + \frac{{b^{2\alpha }  - a^{2\alpha } }}{{2\alpha }}\left( {b - a} \right)} \right\} .
\tag{3.14}
\end{align*}
\paragraph{\textbf{Proof.}}
Multiplying the left hand side and right hand sides of $(3.6)$ and $(3.7)$ we have
\begin{align*}
&t^{2\alpha } f\left( x \right)g\left( x \right) - \left( {t^\alpha  f\left( t \right)} \right)\left( {x^\alpha  g\left( x \right)} \right) \\
&- \left( {x^\alpha  f\left( x \right)} \right)\left( {t^\alpha  g\left( t \right)} \right) + x^{2\alpha } f\left( t \right)g\left( t \right) \\ 
& =\left( {\alpha f\left( \xi  \right) - \xi ^{2 - \alpha } D_\alpha  f\left( \xi  \right)} \right)\left( {\alpha g\left( \eta  \right) - \eta ^{2 - \alpha } D_\alpha  g\left( \eta  \right)} \right)\left( {\frac{{t^\alpha   - x^\alpha  }}{\alpha }} \right)^2 .
\tag{3.15}
\end{align*}
Integrating both sides of $(3.15)$ with respect to $t$ over $[a,b]$ we have
\begin{align*}
&\left( {\frac{{b^{3\alpha }  - a^{3\alpha } }}{{3\alpha }}} \right)f(x)g(x) - x^\alpha  g\left( x \right)\int\limits_a^b {t^\alpha  f\left( t \right)d_\alpha  } t \\
&- x^\alpha  f\left( x \right)\int\limits_a^b {t^\alpha  g\left( t \right)d_\alpha  } t + x^{2\alpha } \int\limits_a^b {f(t)g(t)} d_\alpha  t \\ 
&= \left( {\alpha f\left( \xi  \right) - \xi ^{2 - \alpha } D_\alpha  f\left( \xi  \right)} \right)\left( {\alpha g\left( \eta  \right) - \eta ^{2 - \alpha } D_\alpha  g\left( \eta  \right)} \right) \\ 
&\left\{ {\frac{{b^{3\alpha }  - a^{3\alpha } }}{{3\alpha }} - 2x^\alpha  \frac{{b^{2\alpha }  - a^{2\alpha } }}{{2\alpha }} + x^{2\alpha } \left( {b - a} \right)} \right\}.
\tag{3.16}
\end{align*}
Now integrating $(3.16)$ with respect to $x$ over $[a,b]$ we have

\begin{align*}
&\left( {\frac{{b^{3\alpha }  - a^{3\alpha } }}{{3\alpha }}} \right)\int\limits_a^b {f(x)g(x)} d_\alpha  x - \left( {\int\limits_a^b {x^\alpha  g\left( x \right)d_\alpha  x} } \right)\left( {\int\limits_a^b {t^\alpha  f\left( t \right)d_\alpha  t} } \right) \\
& - \left( {\int\limits_a^b {x^\alpha  f\left( x \right)d_\alpha  x} } \right)\left( {\int\limits_a^b {t^\alpha  g\left( t \right)d_\alpha  t} } \right) + \left( {\int\limits_a^b {x^{2\alpha } d_\alpha  x} } \right)\left( {\int\limits_a^b {f(t)g(t)} d_\alpha  t} \right) \\ 
& = \left( {\alpha f\left( \xi  \right) - \xi ^{2 - \alpha } D_\alpha  f\left( \xi  \right)} \right)\left( {\alpha g\left( \eta  \right) - \eta ^{2 - \alpha } D_\alpha  g\left( \eta  \right)} \right) \\
&\left\{ {\frac{{b^{3\alpha }  - a^{3\alpha } }}{{3\alpha }}\left( {b - a} \right) - \frac{{\left( {b^{2\alpha }  - a^{2\alpha } } \right)^2 }}{{2\alpha }} + \frac{{b^{2\alpha }  - a^{2\alpha } }}{{2\alpha }}\left( {b - a} \right)} \right\}.
\tag{3.17}
\end{align*}
From the properties of modulus and $(3.17)$ we have
\begin{align*}
 \left| {H\left( {f,g} \right)} \right| 
&\le \left| {\alpha f\left( \xi  \right) - \xi ^{2 - \alpha } D_\alpha  f\left( \xi  \right)} \right| \times  \\
& \left| {\alpha g\left( \eta  \right) - \eta ^{2 - \alpha } D_\alpha  g\left( \eta  \right)} \right|\left( {\frac{{3\alpha }}{{b^{3\alpha }  - a^{3\alpha } }}} \right) \times  \\
&\left\{ {\frac{{b^{3\alpha }  - a^{3\alpha } }}{{3\alpha }}\left( {b - a} \right) - \frac{{\left( {b^{2\alpha }  - a^{2\alpha } } \right)^2 }}{{2\alpha ^2 }} + \frac{{b^{2\alpha }  - a^{2\alpha } }}{{2\alpha }}\left( {b - a} \right)} \right\} .
\tag{3.18}
\end{align*}
which is required inequality.
 \paragraph{\textbf{Acknowledgment}}
This research is supported by Science and Engineering Research Board (SERB, New
Delhi, India) Project File No. SB/S4/MS:861/13.
  
\end{document}